\documentclass[a4paper]{amsart}

\pdfoutput=1

\usepackage{amsmath,amssymb,amsthm,url,scalerel}
\usepackage[utf8]{inputenc}
\usepackage[T1]{fontenc}
\usepackage{libertine}
\usepackage[libertine,cmintegrals,cmbraces]{newtxmath}
\usepackage{verbatim,enumitem}
\usepackage{stmaryrd}
\usepackage{mathtools}
\usepackage{microtype}
\usepackage{xcolor}
\definecolor{darkred}{RGB}{139,0,0}
\definecolor{darkblue}{RGB}{0,0,139}
\definecolor{darkgreen}{RGB}{0,100,0}
\usepackage{hyperref}
\hypersetup{
  colorlinks   = true, 
  urlcolor     = darkred, 
  linkcolor    = darkred, 
  citecolor   = darkgreen}
\usepackage{tikz}
\usepackage{tikz-cd}
\usepackage[enableskew]{youngtab}
\usepackage{ytableau}
\usepackage[capitalise]{cleveref}
\usepackage{nicefrac}

\AtBeginDocument{%
   \def\MR#1{}
}
\setenumerate[1]{label=(\roman*),}
\setitemize[1]{label=\raisebox{0.25ex}{\tiny$\bullet$}}

\newcommand{\Diffuo}{\ensuremath{\operatorname{Diff}}}
\newcommand{\Diff}{\ensuremath{\operatorname{Diff}^{\scaleobj{0.8}{+}}}}

\newcommand{\bP}{\operatorname{bP}}
\newcommand{\Sp}{\operatorname{Sp}}
\newcommand{\BSp}{\operatorname{BSp}}
\newcommand{\BU}{\operatorname{BU}}
\newcommand{\sgn}{\operatorname{sgn}}
\newcommand{\Hom}{\operatorname{Hom}}
\newcommand{\coker}{\operatorname{coker}}

\newcommand{\BG}{\operatorname{BG}}

\DeclareMathAlphabet{\mathpzc}{OT1}{pzc}{m}{it}
\newcommand{\catsingle}[1]{\ensuremath{\mathcal{#1}}}

\newcommand{\oH}{\ensuremath{\operatorname{H}}}

\newcommand{\bfR}{\ensuremath{\mathbf{R}}}
\newcommand{\bfZ}{\ensuremath{\mathbf{Z}}}

\newcommand{\cH}{\ensuremath{\catsingle{H}}}

\newcommand{\ra}{\rightarrow}
\newcommand{\lra}{\longrightarrow}

\newcommand{\xra}[1]{\xrightarrow{#1}}
\newcommand{\xlra}[1]{\overset{#1}{\longrightarrow}}

\newcommand{\Mod}[1]{\ (\mathrm{mod}\ #1)}

\newtheorem*{nthm}{Theorem}

\theoremstyle{definition}

\theoremstyle{remark}

\newtheorem*{nrem}{Remark}

\newcommand{\orw}[1]{\marginpar{\tiny\textcolor{blue}{orw: #1}}}
\newcommand{\mk}[1]{\marginpar{\tiny\textcolor{olive}{mk: #1}}}

\begin{document}

\title{Mapping class groups of simply connected\\high-dimensional manifolds need not be arithmetic}

\author{Manuel Krannich}
\email{krannich@dpmms.cam.ac.uk}

\author{Oscar Randal-Williams}
\email{o.randal-williams@dpmms.cam.ac.uk}
\address{Centre for Mathematical Sciences, Wilberforce Road, Cambridge CB3 0WB, UK}


\begin{abstract}
It is well known that Sullivan showed that the mapping class group of a simply connected high-dimensional manifold is commensurable with an arithmetic group, but the meaning of ``commensurable'' in this statement seems to be less well known. We explain why this result fails with the now standard definition of commensurability by exhibiting a manifold whose mapping class group is not residually finite. We do not suggest any problem with Sullivan's result: rather we provide a gloss for it.

\end{abstract}

\maketitle

In the landmark paper \cite{SullivanInf} Sullivan proves a structural result on the groups $\pi_0 \mathrm{Diff}(M)$ of isotopy classes of diffeomorphisms of a simply connected compact manifold $M$ of dimension at least $5$, which has as consequence that these groups are ``commensurable'' with arithmetic groups (Theorem 13.3 of loc.\ cit.). 

Sullivan defines ``commensurability'' of groups (at the bottom of p.\ 307 of loc.\,cit.) as the equivalence relation generated by passing to subgroups of finite index and taking quotients by finite normal subgroups. This differs from the current common usage of this term as the equivalence relation generated only by passing to subgroups of finite index; from now on let us reserve \emph{commensurable} for the latter term, and refer to Sullivan's notion as \emph{commensurable up to finite kernel}.

These two notions differ. Recall that a group is \emph{residually finite} if the intersection of all its finite index subgroups is trivial, or equivalently, if each nontrivial element remains nontrivial in some finite quotient. As arithmetic groups are residually finite \cite[p.\,108]{SerreArith}, and residual finiteness is clearly preserved by passing to finite index sub- and supergroups, a group commensurable to an arithmetic group is in particular residually finite. However, work of Deligne \cite{Deligne} shows that there are central extensions
\[
0\lra \bfZ/n \lra \Gamma \lra \Sp_{2g}(\bfZ)\lra0
\]
for which $\Gamma$ is not residually finite. Such a $\Gamma$ is clearly commensurable up to finite kernel with the arithmetic group $\Sp_{2g}(\bfZ)$, but cannot be commensurable with an arithmetic group.

The aim of this note is to explain how Deligne's example may be imported into manifold theory to provide a mapping class group of a simply connected high-dimensional manifold that is not residually finite, and hence not commensurable with an arithmetic group.

\begin{nthm}
For $n\equiv 5\Mod{8}$ and $g \geq 5$, the mapping class group $\pi_0\Diffuo(\sharp^g S^n \times S^n)$ is not residually finite. 
\end{nthm}

The best known family of groups commensurable up to finite kernel with an arithmetic group but not residually finite is Deligne's, and the best known (to the authors, as we have studied them elsewhere \cite{GRWabelian, krannich}) family of mapping class groups of simply connected high-dimensional manifolds are those of $\sharp^g S^n \times S^n$.  Fortunately these examples are often close to each other; the proof of the Theorem will be to make this precise. 

\begin{nrem}
Sullivan also proves that the groups $\pi_0 \mathrm{hAut}(X)$ of homotopy classes of homotopy equivalences of a simply connected finite CW complex $X$ are commensurable up to finite kernel with an arithmetic group \cite[Theorem 10.3 (i)]{SullivanInf}. However, as explained by Serre \cite[p.\,108]{SerreArith}, in this case further results of Sullivan may be used to show that such groups are commensurable to arithmetic groups in our sense. Indeed, the groups $\pi_0 \mathrm{hAut}(X)$ are residually finite as a consequence of \cite[Theorem 3.2]{SullivanGen}, and among residually finite groups being commensurable up to finite kernel and being commensurable agree: using the second description of residual finiteness mentioned above, one can deduce that for any finite normal subgroup $K \le G$ of a residually finite group $G$, there is a finite group $C$ and a morphism $G\ra C$ such that the induced morphism $ G\ra (G/K)\times C$ is injective and hence exhibits $G$ as a finite index subgroup of $(G/K)\times C$, so $G$ and $G/K$ are commensurable.
\end{nrem}

As a point of terminology, the \emph{finite residual} of a group is the intersection of all its finite index subgroups; a group is residually finite precisely if its finite residual is trivial.

\subsection*{Deligne's extensions}We begin by explaining the special case of Deligne's work \cite{Deligne} mentioned above in a bit more detail (see also \cite{Morris}). The fundamental group of the Lie group $\Sp_{2g}(\bfR)$ can be identified with $\bfZ$, so the pullback of its universal cover to $\Sp_{2g}(\bfZ)$ yields a central extension of the form
\begin{equation}\label{equ:Deligne}
0\lra \bfZ \lra \widetilde{\Sp}_{2g}(\bfZ)\lra \Sp_{2g}(\bfZ)\lra 0.
\end{equation} 
As long as $\Sp_{2g}(\bfZ)$ has the congruence subgroup property, i.e.\,for $g\ge2$, Deligne's work implies that the finite residual of $\widetilde{\Sp}_{2g}(\bfZ)$ agrees with the subgroup $2\cdot\bfZ\subset \widetilde{\Sp}_{2g}(\bfZ) $, so $\widetilde{\Sp}_{2g}(\bfZ)$ is in particular not residually finite. Moreover, this shows that for $n\ge3$ the quotients $\widetilde{\Sp}_{2g}(\bfZ) / (n \cdot \bfZ)$ are not residually finite either, and hence give rise to non-residually finite central extensions of $\Sp_{2g}(\bfZ)$ as asserted in in the introduction. 

\medskip

There is an alternative description of the extension \eqref{equ:Deligne}, more convenient for our purposes, which we now explain.

\subsection*{The signature class}

Let us describe Meyer's \emph{signature class} $\sgn\in\oH^2(\BSp_{2g}(\bfZ))$, following \cite{MeyerThesis}. For $g\ge3$ the group $\Sp_{2g}(\bfZ)$ is perfect, i.e.\,the first homology $\oH_1(\BSp_{2g}(\bfZ))$ vanishes, so by the universal coefficient theorem it suffices to specify a morphism
\begin{equation}\label{equ:signature}
\sgn\colon\oH_2(\BSp_{2g}(\bfZ))\lra\bfZ.
\end{equation}
To do so, recall that, as a consequence of classical work of Thom \cite{Thom}, the second integral homology group of a space $X$ is canonically isomorphic to the group of cobordism classes of continuous maps $f\colon S\ra X$ for oriented closed surfaces $S$. As a result, we may represent a homology class in $\oH_2(\BSp_{2g}(\bfZ))$ by a map $f\colon S\ra X$ of this kind and associate to it a symmetric bilinear form given by the composition 
\[
\oH^1(S;f^*\cH)\otimes \oH^1(S;f^*\cH)\xlra{\smile}\oH^2(S;f^*(\cH\otimes\cH))\xra{\lambda_*}\oH^2(S;\bfZ)\xlra{\langle [S],-\rangle}\bfZ.
\]
Here $\cH$ denotes the local system $\cH\ra \BSp_{2g}(\bfZ)$ induced by the  standard module $\bfZ^{2g}$ of $\Sp_{2g}(\bfZ)$, which comes equipped with a canonical fibrewise symplectic form $\lambda\colon \cH\otimes \cH\ra \bfZ$. In \cite[\S 3]{MeyerThesis} Meyer notes that this form is nondegenerate after passing to the torsion free quotient of $\oH^1(S;f^*\cH)$ (as a consequence of Poincaré duality). Moreover, he shows that the signature of the induced non-degenerate form on the torsion free quotient depends only on the bordism class of $f$ (a a consequence of Lefschetz duality), so taking signatures yields a morphism of the form \eqref{equ:signature}.

The signature class is closely related to the extension \eqref{equ:Deligne}: by construction, this extension is pulled back from the central extension of $\BSp_{2g}(\bfR)\simeq\BU(g)$ that, under the classification of central extensions in terms of second cohomology, corresponds to a generator in $\oH^2(\BSp_{2g}(\bfR))\cong\Hom(\pi_1\Sp_{2g}(\bfR),\bfZ)\cong\bfZ$. This group is also generated by the first Chern class $c_1\in\oH^2(\BSp_{2g}(\bfR))$, so up to possibly rechoosing the identification $\pi_1\Sp_{2g}(\bfR)\cong\bfZ$, the extension \eqref{equ:Deligne} is classified by the pullback $c_1\in\oH^2(\BSp_{2g}(\bfZ))$ along the inclusion $\Sp_{2g}(\bfZ)\subset \Sp_{2g}(\bfR)$. It is a consequence of Meyer's signature formula \cite[\S 4]{MeyerThesis} (see also \cite[p.\,246]{Meyer}) that
\begin{equation}\label{equ:Meyeridentity}
\sgn= -4\cdot c_1\in\oH^2(\BSp_{2g}(\bfZ)).
\end{equation}

\subsection*{The theta subgroup}The mapping class groups we shall consider are closely connected to a certain finite index subgroup $\Sp_{2g}^q(\bfZ)\le\Sp_{2g}(\bfZ)$, sometimes called the \emph{theta group}, which consists of those symplectic matrices that preserve the standard quadratic refinement, i.e.\,the function $q\colon \bfZ^{2g}\ra \bfZ/2$ mapping $\sum_i(a_i\cdot e_i+b_i\cdot f_i)$ to $\sum_{i}a_ib_i$, where $(e_i,f_i)_{1\le i\le g}$ is the usual symplectic basis. By \cite[Lemma 7.5]{GRWabelian}, the signature morphism \eqref{equ:signature} becomes divisible by $8$ when restricted to $\oH_2(\BSp^q_{2g}(\bfZ))$, but the resulting morphism $\sgn/8\colon \oH_2(\BSp^q_{2g}(\bfZ))\ra \bfZ$ does not uniquely determine a class in $\oH^2(\BSp^q_{2g}(\bfZ))$, because in contrast to $\Sp_{2g}(\bfZ)$, the subgroup $\Sp_{2g}^q(\bfZ)$ is not perfect. To specify such a class, we use that the abelianisation of $\Sp_{2g}^q(\bfZ)$ is by \cite[p.\,147]{JohnsonMillson} cyclic of order $4$ and generated by 
\[
X_g=\operatorname{diag}\left(\left(\begin{smallmatrix}0&-1\\1&0\end{smallmatrix}\right),\left(\begin{smallmatrix}1&0\\0&1\end{smallmatrix}\right),\ldots,\left(\begin{smallmatrix}1&0\\0&1\end{smallmatrix}\right)\right)
\]
as long as $g\ge 3$. This implies that the morphism $\bfZ/4\ra \Sp_{2g}^q(\bfZ)$ defined by $X_g$ induces an isomorphism on torsion subgroups of second cohomology, so induces a splitting
\begin{equation}\label{equ:uct}
\oH^2(\BSp_{2g}(\bfZ))\xlra{\cong}\Hom(\oH_2(\BSp_{2g}(\bfZ)),\bfZ)\oplus \oH^2(\mathrm{B}\bfZ/4),
\end{equation}
which we can use to define a class 
\[
\mu\in \oH^2(\BSp^q_{2g}(\bfZ))
\]
by declaring it to agree with the divided signature $\sgn/8$ on the first summand and to be trivial on the second. This class is related to the pullback of $c_1\in \oH^2(\BSp_{2g}(\bfZ))$ by 
\begin{equation}\label{equ:relation}
8\cdot \mu= -4\cdot c_1\in \oH^2(\BSp_{2g}^q(\bfZ)),
\end{equation}
since both sides agree in the first summand of \eqref{equ:uct} as a result of \eqref{equ:Meyeridentity} and in the second because $\oH^2(\mathrm{B}\bfZ/4)$ is annihilated by $4$. By Deligne's result the finite residual of the central extension \eqref{equ:Deligne} classified by $c_1$ is $2\cdot\bfZ$, which implies that the analogous central extension of $\Sp_{2g}(\bfZ)$ classified by $-4 \cdot c_1$ has finite residual $8\cdot\bfZ$, using that
\begin{enumerate}
\item for a central extension $0\ra A\ra E\ra G\ra 0$ classified by a class $\lambda\in \oH^2(\BG;A)$, the multiple $k\cdot \lambda\in \oH^2(\BG;A)$ for $k\in\bfZ$ classifies the pushout of this extension along the multiplication map $k\cdot(-)\colon A\ra A$ (see \cite[p.\,94]{Brown}), and that
\item taking finite residuals is functorial with respect to all group homomorphisms, and maps inclusions of finite index to isomorphisms. 
\end{enumerate}
Since $\Sp_{2g}^q(\bfZ) \le\Sp_{2g}(\bfZ)$ has finite index and we just argued that the central extension of $\Sp_{2g}(\bfZ)$ classified by $-4 \cdot c_1\in\oH^2(\BSp_{2g}(\bfZ))$ has finite residual $8\cdot\bfZ$, we conclude that that the central extension of $\Sp_{2g}^q(\bfZ)$ classified by $-4 \cdot c_1\in\oH^2(\BSp^q(\bfZ))$ also has finite residual $8\cdot\bfZ$. Using (i) and  (ii) once more, we derive from \eqref{equ:relation} that the central extension 
\begin{equation}\label{equ:extSpq}
0\lra \bfZ \lra E_g\lra \Sp_{2g}^q(\bfZ)\lra 0\end{equation}
classified by $\mu$ has finite residual $\bfZ$.

\begin{nrem}
It does not matter for our argument, but the relation \eqref{equ:relation} cannot be divided by 4. By definition $\mu$ vanishes when restricted to $\mathrm{B}\bfZ/4$, but $c_1$ restricts to a generator of $\oH^2(\mathrm{B}\bfZ/4) \cong \bfZ/4$ because the composition $\mathrm{B}\bfZ/4 \to \BSp_{2}^q(\bfZ) \to \BSp_{2}(\bfR) \simeq \BU(1)$ classifies the tautological line bundle.
\end{nrem}

\subsection*{Proof of the Theorem}
Using work of Kreck \cite{Kreck}, it is shown in 
 \cite[Equation 7.3 and Theorem 7.7]{GRWabelian} that for $g\ge5$ and $n\equiv5\Mod{8}$ the mapping class group of orientation preserving diffeomorphisms of $W_g = \sharp^g S^n\times S^n$ admits a splitting 
\[
\pi_0\Diff(W_g)\cong \widetilde{E}_g\times \coker(J)_{2n+1}
\]
into the (finite) cokernel $\coker(J)_{2n+1}$ of the stable $J$-homomorphism in dimension $2n+1$ and the central extension
\[
0\lra \bfZ/\lvert\bP_{2n+2}\rvert\lra \widetilde{E}_g\lra \Sp_{2g}^q(\bfZ)\lra 0,
\]
obtained from \eqref{equ:extSpq} by taking quotients by the subgroup $\lvert\bP_{2n+2}\rvert\cdot \bfZ\le \bfZ$ generated by the order of Kervaire--Milnor's finite cyclic group $\bP_{2n+2}$ of homotopy $(2n+1)$-spheres bounding a parallelisable manifold \cite{KervaireMilnor}. 
Since the finite residual of $E_g$ is $\bfZ$ by the discussion above, that of $\widetilde{E}_g$ is $\bfZ/\lvert\bP_{2n+2}\rvert$, which is nontrivial, since under our assumption that $n\equiv5\Mod{8}$, the group $\bP_{2n+2}$ is known not to vanish (see Section 7 loc.\,cit.). This implies that $ \widetilde{E}_g$ is not residually finite, and hence neither is $\pi_0\Diff(W_g)$ nor its non-orientation preserving variant $\pi_0\Diffuo(W_g)$.

\subsection*{Acknowledgements} The authors were supported by the ERC under the European Union's Horizon 2020 research and innovation programme (grant agreement No.\ 756444), and by a Philip Leverhulme Prize from the Leverhulme Trust.

\bibliographystyle{amsalpha}
\bibliography{literature}

\vspace{0.2cm}

\end{document}